\documentclass[12pt]{neumann}
\usepackage{subfigure}

\author{Benjamin J. Wilson}
\title{Representations of Truncated Current Lie Algebras}

\usepackage{graphicx}

\theoremstyle{definition}

\newtheorem*{theoremnonum}{Theorem}

\newcommand{\UEA}[1]{\text{U}(#1)}
\newcommand{\K}{\Bbbk}

\newcommand{\Z}{\mathbb{Z}}

\newcommand{\NewTerm}[1]{\textit{#1}}
\newcommand{\Indet}[1]{#1} 

\newcommand{\SFrac}[2]{\textstyle{\frac{#1}{#2}}}
\newcommand{\Span}{\text{span}}

\newcommand{\LieBrac}[2]{\lbrack \hspace{0.15em} #1 , #2 \hspace{0.15em}\rbrack} 

\newcommand{\LinearDual}[1]{#1^{*}}
\newcommand{\g}{\mathfrak{g}}
\newcommand{\h}{\mathfrak{h}}
\newcommand{\Roots}{\Delta}
\newcommand{\PositiveRoots}{\Roots_{+}}
\newcommand{\NegativeRoots}{\Roots_{-}}
\newcommand{\RootSpace}[2]{{#1}^{#2}}

\newcommand{\gHat}{\hat{\g}}
\newcommand{\hHat}{\hat{\h}}
\newcommand{\speciallinear}[2]{\mathrm{sl}(#1, #2)}
\newcommand{\Det}[1]{\text{det}\hspace{0.1em} #1}
\newcommand{\Verma}[1]{\mathbf{M}(#1)}
\newcommand{\VermaGen}[1]{\mathbf{v}_{#1}}

\newcommand{\Nilp}{N} 
\newcommand{\hElementMap}{h}
\newcommand{\hElement}[1]{\hElementMap_{#1}} 

\newcommand{\VertForm}[2]{(\hspace{0.2em}#1 \hspace{0.2em} \vert \hspace{0.2em}#2 \hspace{0.2em})}
\newcommand{\Vir}{\mathfrak{a}}
\newcommand{\VirL}{L}
\newcommand{\VirC}{c}

\newcommand{\Functional}[2]{\langle #1, #2 \rangle}
\begin{document}
\begin{abstract}
Let $\g$ denote a Lie algebra, and let $\gHat$ denote the tensor product of $\g$
with a ring of truncated polynomials. The Lie algebra $\gHat$ is called a
truncated current Lie algebra.
The highest-weight theory of $\gHat$ is investigated, and a reducibility criterion for the Verma modules is described.
\end{abstract}
\maketitle
Let $\g$ be a Lie algebra over a field $\K$ of characteristic zero, and fix a positive integer $\Nilp$.
The Lie algebra
\begin{equation}\label{AffinizationEqn}
\gHat = \g \otimes_\K \K[\Indet{t}] / {\Indet{t}^{\Nilp+1} \K[\Indet{t}]},
\end{equation}
over $\K$, with the Lie bracket given by
\begin{equation*}
\LieBrac{x \otimes \Indet{t}^i}{y \otimes \Indet{t}^j} = \LieBrac{x}{y} \otimes \Indet{t}^{i+j} \quad \mbox{for all $x,y \in \g$ and $ i,j \geqslant 0$},
\end{equation*}
is called a \NewTerm{truncated current Lie algebra}, 
or sometimes a \NewTerm{generalised Takiff algebra} or a \NewTerm{polynomial Lie algebra}. 
We describe a highest-weight theory for $\gHat$, and the reducibility criterion for the universal objects of this theory, the Verma modules. 
Representations of truncated current Lie algebras have been studied in \cite{GeoffriauOne, GeoffriauTwo, RaisTauvel, Takiff}, and have applications in the theory of soliton equations \cite{CasatiOrtenzi} and in the representation theory of affine Kac--Moody Lie algebras \cite{WilsonPhDThesis}.

A highest-weight theory is defined by a choice of \NewTerm{triangular decomposition}.
Choose an abelian subalgebra $\h \subset \g$ that acts diagonally upon $\g$ via the adjoint action, and write
\[
\textstyle \g = \h \oplus (\bigoplus_{\alpha \in \Roots}{\RootSpace{\g}{\alpha}})
\]
for the eigenspace decomposition, where $\Roots \subset \LinearDual \h$, and for all $\alpha \in \Roots$,
\[
\LieBrac{h}{x} = \Functional{\alpha}{h} \hspace{0.1em} x, \quad \mbox{for all $h \in \h$ and $x \in \RootSpace\g\alpha$}.
\]
A triangular decomposition of $\g$ is, in essence\footnote{There are additional hypotheses (see \cite{MoodyPianzola}) --- in particular, there must exist some finite subset of $\PositiveRoots$ that generates $\PositiveRoots$ under addition.
This excludes, for example, the imaginary highest-weight theory of an affine Lie algebra (see \cite{WilsonTCLA}, Appendix B).}, a division of the eigenvalue set $\Roots$ into two opposing halves
\begin{equation}\label{RootSpaceDecomposition}
\Roots = \PositiveRoots \sqcup \NegativeRoots, \qquad -\PositiveRoots = \NegativeRoots,
\end{equation}
that are closed under addition, in the sense that the sum of two elements of $\PositiveRoots$ is another element of $\PositiveRoots$, if it belongs to $\Roots$ at all. 
The decomposition \eqref{RootSpaceDecomposition} defines a decomposition of $\g$ as a direct sum of subalgebras
\begin{equation}\label{TriangularDecomposition}
\textstyle \g = \g_{+} \oplus \h \oplus \g_{-}, \quad \text{where} \quad \g_{\pm} = \bigoplus_{\alpha \in \PositiveRoots}{\g^{\pm \alpha}}.
\end{equation}
For example, if $\g = \speciallinear 3 \K$, the Lie algebra of traceless $3 \times 3$ matrices with entries from the field $\K$,
then the subalgebras $\h$, $\g_{+}$, $\g_{-}$ are the traceless diagonal, upper-triangular and lower-triangular matrices, respectively.
In analogy with the classical case, where $\g$ is finite-dimensional and semisimple, $\h$ might be called a \NewTerm{diagonal subalgebra} or a \NewTerm{Cartan subalgebra}, while the elements of $\Roots$ and $\PositiveRoots$ might be called \NewTerm{roots} and \NewTerm{positive roots}, respectively.

The concept of a triangular decomposition is also applicable to many infinite-dimensional Lie algebras of importance in mathematical physics, such as Kac--Moody Lie algebras, the Virasoro algebra and the Heisenberg algebra.
For example, the \NewTerm{Virasoro algebra} is the $\K$-vector space $\Vir$ with basis the set of symbols 
\[
\set{\VirL_m | m \in \Z} \cup \set{ \VirC},
\]
endowed with the Lie bracket given by
\[
\LieBrac{\VirC}{\Vir}=\set{0}, \qquad \LieBrac{\VirL_m}{\VirL_n} = (m-n) \VirL_{m+n} + \delta_{m,-n} \SFrac{m^3 - m}{12} \VirC, 
\]
for all $m, n \in \Z$, where $\delta$ denotes the Kronecker function.
If $\g = \Vir$, then the subalgebras
\[
\h = {\K \VirL_0} \oplus {\K \VirC}, \qquad \g_{\pm} = \Span \set{ \VirL_{\pm m} | m > 0},
\]
provide a triangular decomposition.

The triangular decomposition \eqref{TriangularDecomposition} of $\g$ naturally defines a triangular decomposition of $\gHat$,
\[ 
\textstyle \gHat = \gHat_{-} \oplus \hHat \oplus \gHat_{+}, \quad \mbox{where} \quad \gHat_{\pm} = \bigoplus_{\alpha \in \PositiveRoots }{\gHat^{\pm \alpha}},
\]
where the subalgebra $\hHat$ and the subspaces $\RootSpace\gHat\alpha$ are defined in the manner of \eqref{AffinizationEqn}, and $\h \subset \hHat$ is the diagonal subalgebra.
Hence a $\gHat$-module $M$ is a \NewTerm{weight module} if the action of $\h$ on $M$ is diagonalisable.
A weight $\gHat$-module is of \NewTerm{highest weight} if there exists a non-zero vector $v \in M$, and a functional $\Lambda \in \LinearDual{\hHat}$ such that
\[ 
\gHat_{+} \cdot v = 0; \qquad \UEA\gHat \cdot v = M; \qquad h \cdot v = \Lambda (h) v \quad \mbox{for all $\quad h \in \hHat$}. 
\]
The unique functional $\Lambda \in \LinearDual{\h}$ is the \NewTerm{highest weight} of the highest-weight module $M$. 
Notice that the weight lattice of a weight module is a subset of $\LinearDual\h$, while a highest-weight is an element of $\LinearDual{\hHat}$.
A highest-weight $\Lambda \in \LinearDual{\hHat}$ may be thought of as a tuple of functionals on $\h$,
\begin{equation}\label{Tuples}
\Lambda = (\Lambda_0, \Lambda_1, \dots, \Lambda_\Nilp) \quad \text{where} \quad \Functional {\Lambda_i} {h} = \Functional {\Lambda} {h \otimes \Indet{t}^i} \quad \mbox{for all $h \in \h$ and $i \geqslant 0$}.
\end{equation}
All $\gHat$-modules of highest weight $\Lambda \in \LinearDual{\hHat}$ are homomorphic images of a certain universal $\gHat$-module of highest weight $\Lambda$, denoted by $\Verma\Lambda$.
These universal modules $\Verma\Lambda$ are known as \NewTerm{Verma modules}.

A single hypothesis suffices for the derivation of a criterion for the reducibility of a Verma module $\Verma\Lambda$ for $\gHat$ in terms of the functional $\Lambda \in \LinearDual{\hHat}$.
We assume that the triangular decomposition of $\g$ is \NewTerm{non-degenerately paired}, i.e., that for each $\alpha \in \PositiveRoots$, a non-degenerate bilinear form
\[ {\VertForm\cdot\cdot}_\alpha : \RootSpace\g\alpha \times \RootSpace\g{-\alpha} \to \K,
\]
and a non-zero element $\hElement\alpha \in \h$ are given, such that
\[ 
\LieBrac{x}{y} = {\VertForm{x}{y}}_\alpha \hspace{0.15em} \hElement\alpha,
\]
for all $x \in \RootSpace\g\alpha$ and $y \in \RootSpace\g{-\alpha}$.
All the examples of triangular decompositions considered above satisfy this hypothesis.
The reducibility criterion is given by the following theorem, which we state without proof.\par
\begin{figure}
  \begin{center}
    \subfigure[$\g = \speciallinear 3 \K$]{\label{A2Figure}\includegraphics[scale=0.9]{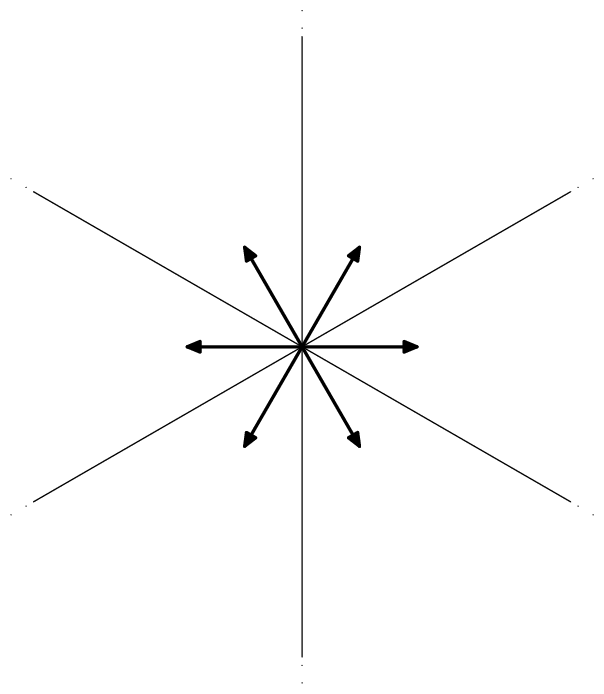}}
    \subfigure[$\g$ the Virasoro algebra $\Vir$]{\label{VirasoroFigure}\includegraphics[scale=0.9]{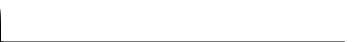}} 
  \end{center}
  \caption{Reducibility criterion for the Verma modules of $\gHat$}
\end{figure}

\vspace{0.5cm}
\begin{theoremnonum}\cite{WilsonTCLA}
The Verma module $\Verma\Lambda$ for $\gHat$ is reducible if and only if
\[
\Functional\Lambda{\hElement\alpha \otimes \Indet{t}^\Nilp} = 0
\]
for some positive root $\alpha \in \PositiveRoots$ of $\g$.
\end{theoremnonum}
\vspace{0.5cm}
Notice that the reducibility of $\Verma\Lambda$ depends only upon $\Lambda_\Nilp \in \LinearDual\h$, the last component of the tuple \eqref{Tuples}.
The criterion described by the theorem has many disguises, depending upon the underlying Lie algebra $\g$.
If $\g = \speciallinear 3 \K $, then $\Verma\Lambda$ is reducible if and only if $\Lambda_\Nilp$ is orthogonal to a root.
This is precisely when $\Lambda_\Nilp$ belongs to one of the three hyperplanes in $\LinearDual\h$ illustrated in Figure \ref{A2Figure}.
The arrows describe the root system.
If instead $\g$ is the Virasoro algebra $\Vir$, then $\Verma\Lambda$ is reducible if and only if
\[
2m \Functional{\Lambda_\Nilp}{\VirL_0} + \SFrac{m^3 - m}{12} \Functional{\Lambda_\Nilp}{\VirC} = 0,
\]
for some non-zero integer $m$. 
That is, $\Verma{\Lambda}$ is reducible when $\Lambda_\Nilp$ belongs to the infinite union of hyperplanes indicated in Figure \ref{VirasoroFigure}.
The extension of a functional in the horizontal and vertical directions is determined by evaluations at $\VirC$ and $\VirL_0$, respectively.

It is remarkable that there exists a unified reducibility criterion for the Verma modules for truncated current Lie algebras $\gHat$, irrespective of the particular Lie algebra $\g$.
The technology available for deriving reducibility criteria, called the Shapovalov determinant, is easier to operate for truncated current Lie algebras.
To illustrate this, we construct the Shapovalov determinant in the special case where $\g$ is finite-dimensional and semisimple.
A Verma module is generated by the action of the subalgebra $\g_{-}$ upon the highest-weight vector.
We might visualise the weight lattice of a Verma module as a cone extending downwards from the highest weight $\Lambda$.
The Poincar\'e--Birkhoff--Witt Theorem informs us that $\Verma\Lambda$ has a basis parameterised by the (unordered) downward paths in the weight lattice beginning at the highest weight; alternatively, these might be conceived of as multisets with entries from $\NegativeRoots$.
If we fix some weight $\Lambda - \chi$, there are only finitely many downward paths, say $\theta_1, \theta_2, \dots, \theta_l$, extending from $\Lambda$ to $\Lambda - \chi$.
Dual to each downward path $\theta_i$ is the upward path $\theta^i$ from $\Lambda - \chi$ to $\Lambda$, obtained by reversing the direction of the arrows.
The highest-weight space of $\Verma\Lambda$ is one-dimensional, spanned by the highest-weight vector $\VermaGen\Lambda$. 
Hence descending from $\VermaGen\Lambda$ via $\theta_i$, and then ascending again via $\theta^j$, results in a scalar multiple $\theta_i^j$ of $\VermaGen\Lambda$.
We consider all such values together, as a matrix $T_\chi^\Lambda$:
\[
T_\chi^\Lambda = (\theta_i^j)_{1 \leqslant i, j \leqslant l} \quad \mbox{where} \quad 
(\theta^j \circ \theta_i) \VermaGen\Lambda = \theta_i^j \VermaGen\Lambda.
\]
The determinant $\Det{T_\chi^\Lambda}$ of this matrix is the \NewTerm{Shapovalov determinant} of $\Verma\Lambda$ at $\chi$.
Degeneracy of this matrix would indicate the existence of a non-zero element of weight $\Lambda - \chi$ that vanishes along all paths to the highest-weight space; such vectors generate proper submodules.
In particular, the Verma module $\Verma\Lambda$ is reducible if and only if $\Det{T_\chi^\Lambda} = 0$ for some $\chi$. 

In the case of the truncated current Lie algebras, there is a straightforward and unified approach to the derivation of a formula for the Shapovalov determinant.
The nilpotency of the indeterminate $\Indet{t}$ is such that many of the would-be non-zero entries of the Shapovalov matrix vanish.
This permits the diagonalisation of the Shapovalov matrix by making a clever choice for the basis $\theta_i$, and redefining the duality between the downward paths $\theta_i$ and the upward paths $\theta^j$.
The calculation of the determinant of such a matrix is only as difficult as the calculation of its diagonal entries.

\bibliographystyle{abbrv}
\bibliography{thesis}

\begin{thebibliography}{1}

\bibitem{CasatiOrtenzi}
P.~Casati and G.~Ortenzi.
\newblock New integrable hierarchies from vertex operator representations of
  polynomial {L}ie algebras.
\newblock {\em J. Geom. Phys.}, 56(3):418--449, 2006.

\bibitem{GeoffriauOne}
F.~Geoffriau.
\newblock Sur le centre de l'alg\`ebre enveloppante d'une alg\`ebre de
  {T}akiff.
\newblock {\em Ann. Math. Blaise Pascal}, 1(2):15--31 (1995), 1994.

\bibitem{GeoffriauTwo}
F.~Geoffriau.
\newblock Homomorphisme de {H}arish-{C}handra pour les alg\`ebres de {T}akiff
  g\'en\'eralis\'ees.
\newblock {\em J. Algebra}, 171(2):444--456, 1995.

\bibitem{MoodyPianzola}
R.~V. Moody and A.~Pianzola.
\newblock {\em Lie algebras with triangular decompositions}.
\newblock Canadian Mathematical Society Series of Monographs and Advanced
  Texts. John Wiley \& Sons Inc., New York, 1995.

\bibitem{RaisTauvel}
M.~Ra{\"{\i}}s and P.~Tauvel.
\newblock Indice et polyn\^omes invariants pour certaines alg\`ebres de {L}ie.
\newblock {\em J. Reine Angew. Math.}, 425:123--140, 1992.

\bibitem{Takiff}
S.~J. Takiff.
\newblock Rings of invariant polynomials for a class of {L}ie algebras.
\newblock {\em Trans. Amer. Math. Soc.}, 160:249--262, 1971.

\bibitem{WilsonTCLA}
B.~J. Wilson.
\newblock Highest-weight theory for truncated current {L}ie algebras.
\newblock arXiv:0705.1203.

\bibitem{WilsonPhDThesis}
B.~J. Wilson.
\newblock {\em Representations of infinite-dimensional {L}ie algebras}.
\newblock PhD thesis, University of Sydney / Universidade de S\~ao Paulo, 2007.
\newblock Submitted.

\end{thebibliography}
\end{document}